\documentclass[12pt]{article}
\usepackage{latexsym,amsmath,amssymb}
\begin{document}
\noindent{\bf \large On the Moment Determinacy of Products of\\
Non-identically Distributed Random Variables}

\vspace{0.8cm} \noindent {\bf Gwo Dong Lin$^{1}$ \ $\bullet$ \
Jordan Stoyanov$^{2}$ }

\vspace{0.2cm}\noindent $^{1}$ Institute of Statistical Science,
Academia Sinica, Taipei 11529, Taiwan (ROC) \\e-mail:
gdlin@stat.sinica.edu.tw

\vspace{0.2cm}\noindent $^{2}$ School of Mathematics $\&$
Statistics, Newcastle University, Newcastle upon Tyne NE1 7RU,
UK\\
e-mail: stoyanovj@gmail.com

\vspace{1.0cm}\noindent {\bf Abstract:} We show first that there
are intrinsic relationships among different conditions, old and
recent, which lead to some general statements in both the
Stieltjes and the Hamburger moment problems. Then we describe
checkable conditions and prove new results about the moment
(in)determinacy for products of independent and non-identically
distributed random variables. We treat all three cases when the
random variables are nonnegative (Stieltjes case), when they take
values in the whole real line (Hamburger case), and the mixed
case. As an illustration we characterize the moment determinacy of
products of random variables whose distributions are generalized
gamma or double generalized gamma all with distinct shape
parameters. Among other corollaries,  the product of two
independent random variables, one  exponential and one inverse
Gaussian, is moment determinate, while the product is moment
indeterminate for the cases: one exponential and one normal, one
chi-square and one normal, and one inverse Gaussian and one
normal.

\vspace{0.2cm}\noindent {\bf Mathematics Subject Classification
(2010)} \ 60E05, 44A60

\vspace{0.2cm}\noindent {\bf Keywords}  Product of random
variables, Stieltjes  moment problem, Hamburger moment problem,
Cram\'{e}r's condition, Carleman's condition, Krein's condition,
Hardy's condition.\\

\vspace{0.4cm}\noindent {\bf 1. Introduction}

\vspace{0.1cm} There is a long standing interest in studying
products of random variables; see, e.g., \cite{l1959},
\cite{l1967}, \cite{s1979}, \cite{cs1994}, \cite{s2000},
\cite{gs2004}, \cite{t2008} and the references therein. The
reasons are two-fold. On one hand, to deal with products leads to
non-trivial, difficult and challenging theoretical problems
requiring to use diverse ideas and techniques. Let us mention just
a few sources: \cite{gs2004}, \cite{b2005}, \cite{b2013}. On the
other hand, products of random variables are naturally involved in
stochastic modelling of complex random phenomena in areas such as
statistical physics, quantum theory, communication theory and
financial modelling; see, e.g., \cite{cm1993}, \cite{fs1997},
\cite{gs2004}, \cite{dm2009}, \cite{hetal2010}, \cite{de2010},
\cite{petal2010}, \cite{cetal2012}.

In general, it is difficult to find explicit
closed-form expressions for the densities or the
distributions of products of random variables with different
distributions. It is, however, possible to
study successfully the moment problem for products of independent
random variables; see, e.g., \cite{ls2014}, \cite{setal2014}. Answers
about the moment (in)determinacy can be found if
requiring only information about the asymptotic
of the moments or about the tails of the densities or of their
distributions.

 All random variables considered in this paper are
defined on an underlying probability space $(\Omega, {\cal F},
{\bf P})$ and we denote by ${\bf E}[X]$ the expected value of the
random variable $X$. A basic assumption is that the random
variables we deal with have finite moments of all positive orders,
i.e. ${\bf E}[|X|^k] < \infty, \ k=1,2,\ldots.$ We write $X\sim F$
to mean that $X$ is a random variable whose distribution function
is $F$ and denote its $k$th order moment by $m_k={\bf E}[X^k].$ We
say that $X$ or $F$ is moment determinate (M-det) if $F$ is the
only distribution having the moment sequence
$\{m_k\}_{k=1}^{\infty}$; otherwise, we say that $X$ or $F$ is
moment indeterminate (M-indet). We use traditional notions,
notations and terms such as Cram\'{e}r's condition, Carleman's
condition, Krein's condition, and Hardy's condition. For reader's
convenience we give their definitions in the text.

We use $\Gamma (\cdot)$ for the Euler-gamma function, ${\mathbb
R}=(-\infty, \infty)$ for the set of all real numbers, ${\mathbb
R}_+=[0,\infty)$ for the nonnegative numbers,  the symbol ${\cal
O}(\cdot)$ with its usual meaning in asymptotic analysis and the
abbreviation i.i.d.\!\! for independent and identically
distributed (random variables).

 In Section 2 we describe useful intrinsic relationships
 among different old and recent conditions
 involved in the Stieltjes and/or the
Hamburger moment problems. Then we prove some new results under
conditions which are relatively easy to check. In Section 3 we
deal with the moment determinacy of products of independent
nonnegative random variables with different distributions, while
in Section 4 we consider products of random variables with values
in ${\mathbb R}$. Finally, in Section 5, we treat the mixed case:
products of both types of random variables, nonnegative ones and
real ones, the latter with values in ${\mathbb R}$.

The results presented in this paper extend some previous results
for products of i.i.d.\! random variables. Here we need a more
refined analysis of the densities of  products than in the
i.i.d.\! case. As an illustration we characterize the moment
(in)determinacy of products of random variables whose
distributions are generalized gamma or double generalized gamma
all with distinct shape parameters. We have derived several
corollaries involving popular distributions widely used in
theoretical studies and applications. Let us list a few: (i) the
product of two independent random variables, one exponential and
one inverse Gaussian, is M-det; (ii) the product of independent
exponential and normal random variables is M-indet; (iii) the
product of independent chi-square and normal random variables is
M-indet; and (iv) the product of independent inverse Gaussian and
normal random variables is M-indet.

\vspace{0.4cm}\noindent
{\bf 2. Some General Results}

\vspace{0.2cm} In this section we present two lemmas,  each
containing workable conditions which, more or less, are available
in the literature. Some of these conditions are old, while others
are recent. We describe intrinsic relationships among these
conditions and use them to prove new results; see Theorems 1--4.

Our findings in this section can be considered as a useful complement to
 the classical criteria of
Cram\'er, Carleman, Krein, Hardy and their converses, so that all
these taken together make more clear, and possibly complete, the
picture of what is in our hands when discussing  the determinacy
of  distributions in terms of their moments.

\vspace{0.2cm}\noindent
{\it 2.1. Stieltjes Case}

\vspace{0.2cm}\noindent
{\bf Lemma 1.} {\it Let $0\le X\sim F$. Then the following four statements
are equivalent:\\
(i) $m_k={\cal O}(k^{2k})$ as $k\to\infty.$ \\
(ii) $\limsup_{k\to\infty}\frac{1}{k}\,m_k^{1/(2k)}<\infty$.\\
(iii) $m_k\le c_0^k\,(2k)!, \ k=1,2,\ldots,$ for some constant $c_0>0$.\\
(iv) $X$ satisfies Hardy's condition, namely, ${\bf
E}[e^{c\sqrt{X}}]<\infty$ for some constant $c>0$.}

\vspace{0.2cm} The equivalence of conditions (i) and (ii), a known
fact for decades, can be easily checked. Conditions (iii) and (iv)
appeared recently and their equivalence to condition (ii) were
shown in \cite{sl2012}.

\vspace{0.2cm}\noindent
{\bf Theorem 1.} {\it Let $0\le X\sim F$ with moments growing as follows:  $m_k={\cal
O}(k^{ak})$ as $k\to\infty$ for some constant $a\in (0, 2]$.
Then the following two statements hold:\\
(i) $X$ satisfies Hardy's condition and hence $X$ is M-det.\\
(ii) The boundary value  $a=2$ is the best possible for $X$ to be
M-det. In fact, there is an M-indet random variable ${\tilde X}\ge
0$ such that ${\bf E}[{\tilde X}^k]={\cal O}(k^{ak})$ as
$k\to\infty$ for all $a>2$.}

\vspace{0.1cm}\noindent {\bf Proof.} Part (i) follows easily from
Lemma 1. To prove part (ii),  we consider the following absolutely
continuous distribution $\tilde F$ whose density is
\begin{eqnarray}
{\tilde f}(x)={\tilde c}\,e^{-\sqrt{x}/(1+|\ln x|^{\delta})},\ \
x>0.
\end{eqnarray}
Here $\delta>1$ and $\tilde c$ is the norming constant. Then it is
easy to evaluate the Krein quantity for ${\tilde X} \sim {\tilde
F}.$ Recall that $\tilde X$ is nonnegative, so in this Stieltjes
case we obtain
$$\hbox{K}[{\tilde f}]:=\int_0^{\infty}\frac{-\ln {\tilde f}
(x^2)}{1+x^2}dx<\infty.$$
Hence $\tilde X$ is M-indet (see, e.g., \cite{l1997}, Theorem 3).

The next step is to check that ${\bf E}[{\tilde X}^k]={\cal
O}(k^{ak})$ as $k\to\infty$ for all $a>2$. To see this, we fix
$a>2,$ take $b\in (2,a)$ and easily find a number $x_0 \geq 1$ such
that $\sqrt{x}>x^{1/b}(1+|\ln x|^{\delta})$ for all $x \ge x_0$. We
now have that
$$\int_0^{x_0}x^k{\tilde f}(x)dx\le \frac{{\tilde c}}{k+1}\,x_0^{k+1}={\cal
O}(k^{ak})~\hbox{ as }~k\to\infty$$ and that
\begin{eqnarray*}
\int_{x_0}^{\infty}x^k{\tilde f}(x)dx&\le&
{\tilde c}\int_{x_0}^{\infty}x^ke^{-x^{1/b}}dx\le
{\tilde c}\int_{0}^{\infty}x^ke^{-x^{1/b}}dx\\&=&b\,{\tilde c}\,\Gamma((k+1)b)={\cal
O}(k^{ak}) \ \hbox{ as } \ k\to\infty.
\end{eqnarray*}
For the last relation we have used the approximation of the gamma function:
\begin{eqnarray*}
\Gamma(x)\approx \sqrt{2\pi}\,x^{x-1/2}\,e^{-x}~~\hbox{ as }~x \rightarrow \infty
\end{eqnarray*}
(see, e.g., \cite{ww1927}, p.\,253). Thus we have shown that indeed
${\bf E}[{\tilde X}^k]={\cal O}(k^{ak})$ as $k\to\infty$ for all
$a>2.$ Therefore the constant 2 in the formulation of the theorem is
indeed the best possible for $X$ to be M-det.   \hspace{\fill}
$\Box$

\vspace{0.15cm}\noindent {\bf Remark 1.}  For $0\le X\sim F$, let
us compare the following two moment conditions: (a) $m_k={\cal
O}(k^{2k})$ as $k\to\infty,$ and (b) $m_{k+1}/m_k={\cal
O}((k+1)^2)$ as $k\to\infty$. Here (a) is the condition in Theorem
1, while condition (b) was introduced and used in the recent paper
\cite{ls2014}. Both conditions are checkable and each of them
guarantees the moment determinacy of $F$. Just to mention that
condition (b) implies condition (a) by referring to Theorem 3 in
\cite{ls2014}, while the converse may not be true in general.

\vspace{0.1cm}
The next result, Theorem 2 below, is the converse to Theorem 1, and deals
with the moment indeterminacy of nonnegative random variables. We need first one condition
which is used a few times in the sequel.

\vspace{0.2cm}\noindent {\bf Condition L:} Suppose, in the
Stieltjes case,  that $f(x), \ x \in {\mathbb R}_+,$ is a density
function such that for some fixed $x_0>0$, \ $f$ is strictly
positive and differentiable for $x > x_0$ and
\begin{eqnarray*}
L_{f}(x):=-\frac{xf^{\prime}(x)}{f(x)}\nearrow
\infty~~\hbox{as}~~x_0<x\rightarrow \infty.
\end{eqnarray*}
In the Hamburger case we require the density
$f(x), \ x \in {\mathbb R},$ to be symmetric.

\vspace{0.1cm} This condition plays a significant r\^ole in moment
problems for absolutely continuous probability distributions.  It
was explicitly introduced and efficiently used for the first time
in \cite{l1997} and later used by several authors naming it as
`Lin's condition'. This condition is involved in some of our
results to follow.

\vspace{0.2cm}\noindent
{\bf Theorem 2.} {\it Let $0\le X\sim F$ and its moment sequence
$\{m_k, \ k=1,2,\ldots \}$  grow `fast' in the sense that
 $m_k\ge
c\,k^{(2+\varepsilon)k}, \ k=1,2,\ldots$, for some constants $c>0$
and $\varepsilon >0$. Assume further that $X$ has a density
function $f$ which satisfies the above Condition L. Then $X$ is
M-indet.}

\vspace{0.2cm}\noindent
{\bf Proof.}  By the condition on the moments, we see that the
Carleman quantity for the moments of $F$ is finite. Indeed,
in this Stieltjes case we have:
$$\mbox{C}[F]= \sum_{k=1}^{\infty}\frac{1}{m_{k}^{1/(2k)}} \le
\sum_{k=1}^{\infty}\frac{1}{c^{1/(2k)}\,k^{1+\varepsilon/2}}<\infty.$$
However no conclusion can be drawn  from this because
$\hbox{C}[F]<\infty$ is only a necessary condition for $X$ to be
M-indet. We need other arguments. By applying Condition L and the
second part of the proof of Theorem 4 in \cite{ls2014}, we finally
conclude that indeed $X$ is M-indet. \hspace{\fill}$\Box$

\vspace{0.15cm}\noindent {\bf Remark 2.} To provide one
application of Theorem 2, let us consider, for example, the random
variable $X=\xi^{2(1+\varepsilon)}$, where $\varepsilon>0$ and
$\xi\sim Exp(1)$, the standard exponential distribution. On one
hand, we can use the Krein criterion and show that $X$ is M-indet.
On the other hand, $X$ satisfies the moment condition in Theorem
2. And here is the point: instead of applying Krein's condition,
we can prove the moment indeterminacy of $X$ by checking that its
density $f$ satisfies Condition L. In general, we follow that
approach which seems easier, or which is working in the specific
case of interest.

\vspace{0.2cm}\noindent
{\it 2.2. Hamburger Case}

\vspace{0.1cm} We start with Lemma 2 establishing the equivalence
of different type of conditions involved to decide whether a
distribution on the whole real line $\mathbb R$ is M-det. Then we
present some new results. Theorem 3 below is a slight
modification, in a new light, of a result in \cite{bs2000},
p.\,92, while Theorem 4 is the converse to Theorem 3.

\vspace{0.2cm}\noindent
{\bf Lemma 2.} {\it Let $X$ be a random variable taking values in
${\mathbb R}$. Then the following
four statements are equivalent:\\
(i) $m_{2k}={\cal O}((2k)^{2k})$ as $k\to\infty$. \\
(ii) $\limsup_{k\to\infty}\frac{1}{2k}\,m_{2k}^{1/(2k)}<\infty$.\\
(iii) $m_{2k}\le c_0^k\,(2k)!,\ k=1,2,\ldots,$ for some constant
$c_0>0$.\\
(iv) $X$ satisfies Cram\'er's condition: its moment generating function exists. }

\vspace{0.2cm}\noindent {\bf Proof.} It is easy to check the
equivalence of conditions (i) and (ii). The equivalence of
conditions (ii) and (iv) is well-known, but we provide here a
simple and instructive proof based on condition (i). Indeed, by
Lemma 1 above, condition (i) is equivalent to say that the random
variable $Y=X^2$ satisfies Hardy's condition, namely, ${\bf
E}[e^{c\sqrt{Y}}]={\bf E}[e^{c|X|}]<\infty$ for some constant
$c>0$. The latter, however,  means that $X$ itself has a moment
generating function. This is exactly the statement (iv). Finally,
applying again Lemma 1 to the nonnegative random variable $Y$, we
obtain the equivalence of (ii) and (iii). Therefore, as stated,
all four conditions (i) -- (iv) are equivalent. \hspace{\fill}
$\Box$

\vspace{0.2cm}\noindent {\bf Remark 3.} In
\cite{bs2000} the moment condition $m_k={\cal O}(k^{k})$ as $k\to \infty$
was used to derive the M-det of $X$  on ${\mathbb R}.$
This condition can be replaced by a weaker one, allowing a `faster' growth of the moments., e.g.,
$m_k={\cal O}((a(k))^k)$ as $k\to \infty,$ where
\[
a(k)=k\,\ln k, \ \mbox{ or } \ k\,(\ln k)\,(\ln \ln k), \mbox{ or
} \ k\,(\ln k)\,(\ln \ln k)\,(\ln \ln \ln k),\ \ldots ,
\]
and still preserving the M-det property of X. Such a statement can
be established by using quasi-analytic functions. We do not give
details here.

\vspace{0.2cm}\noindent
{\bf Theorem 3.} {\it Let $X\sim F,$ where
 $F$ has an unbounded support $\mbox{supp}(F) \subset {\mathbb
R}$ and its moments satisfy the condition: \
$m_{2k}={\cal O}((2k)^{2ak})$ as $k \to \infty$ for some constant
$a\in (0, 1]$.
Then the following statements hold:\\
(i) $X$ satisfies Cram\'er's condition  and hence is M-det. \\
(ii) The boundary value $a=1$ is the best possible for $X$ to be
M-det. In fact, there is an M-indet random variable $ \tilde X$ such
that ${\bf E}[{\tilde X}^{2k}]={\cal O}((2k)^{2ak})$ as $k\to\infty$
for all $a>1$.}

\vspace{0.2cm}\noindent {\bf Proof.} Part (i) follows immediately
from Lemma 2. Part (ii) is essentially given in \cite{bs2000}, but
we provide here a somewhat general case. Let us consider the
distribution $\tilde F$ with the following symmetric density
(compare this with (1)):
\begin{eqnarray}
{\tilde f}(x)={\bar c}\,e^{-|x|/(1+|\ln |x||^{\delta})},\ \ x\in
{\mathbb R},
\end{eqnarray}
where $\delta>1$, ${\tilde f}(0)=1$ and $\bar c$ is a norming
constant. In this, already Hamburger case, the Krein quantity for
${\tilde X}\sim {\tilde F},$ with density $\tilde f$ given by (2),
can be evaluated and shown to be finite. Namely, we have that
$$\hbox{K}[{\tilde f}]:=\int_{-\infty}^{\infty}\frac{-\ln {\tilde f}(x)}{1+x^2}dx<\infty.$$
Hence $\tilde X$ is M-indet (see, e.g., \cite{l1997}, Theorem 1).
However, it is seen that ${\bf E}[{\tilde X}^{2k}]={\cal
O}((2k)^{2ak})$ as $k\to\infty$ for all $a>1$. Therefore, the
constant $1$ in the formulation of the theorem is indeed the best
possible for $X$ to be M-det. \hspace{\fill} $\Box$

\vspace{0.1cm}\noindent {\bf Remark 4.}  Suppose  $X\sim F$
with $F$ having unbounded support, $\mbox{supp}(F) \subset {\mathbb
R}$. We want to compare the following two moment conditions: (a)
$m_{2k}={\cal O}((2k)^{2k})$ as $k\to\infty$, and (b)
$m_{2(k+1)}/m_{2k}={\cal O}((k+1)^2)$ as $k\to\infty$. Here (a) is
the growth of the moments condition in  Theorem 3, while condition
(b) was introduced and successfully exploited in the recent work
\cite{setal2014}. Both conditions are checkable and each of them
guarantees the moment determinacy of $X$ and $F$. Let us mention
that condition (b) implies condition (a) by referring to Theorem 2
in \cite{setal2014}, while the converse may not in general be
true.

\vspace{0.2cm}\noindent {\bf Theorem 4.} {\it Suppose that the
moments of $X\sim F$ grow `fast' in the sense that $m_{2k}\ge
c(2k)^{2(1+\varepsilon)k}, \ k=1,2,\ldots,$ for some positive
constants $c$ and $\varepsilon$. Assume further that $X$ has a
density function $f$ which is symmetric about zero and  satisfies
the above Condition L. Then $X$ is M-indet.}

\vspace{0.2cm}\noindent
{\bf Proof.} By the condition on the moments, we see that
the Carleman quantity for  $F$ is finite
(we remember that this is a Hamburger case):
$$\hbox{C}[F]:= \sum_{k=1}^{\infty}\frac{1}{m_{2k}^{1/(2k)}} \le
\sum_{k=1}^{\infty}\frac{1}{c^{1/(2k)}\,(2k)^{1+\varepsilon}}<\infty.$$
Since no conclusion can be drawn from this finding, we need
different arguments. We use Condition L and the second part of the
proof of Theorem 3 in \cite{setal2014} thus arriving at the
desired conclusion that indeed $X$ is M-indet.\hspace{\fill}
$\Box$

\vspace{0.2cm}\noindent {\bf Remark 5.} For example, instead of
applying Krein's condition, we can use Theorem 4 to prove the
moment indeterminacy of $X\sim F$ whose density is the
symmetrization of that of $\xi^{1+\varepsilon}$, where
$\varepsilon>0$ and $\xi\sim Exp(1)$.

\vspace{0.4cm}\noindent {\bf 3. Products of Nonnegative Random Variables}

\vspace{0.2cm}
We start with two results describing relatively simple
conditions on the
random variables $\xi_1, \ldots, \xi_n$ in order to
guarantee that their product is M-det.

\vspace{0.2cm}\noindent
{\bf Theorem 5.} {\it Suppose that the moments
$m_{i,k}={\bf E}[\xi_i^k], \ i=1,\ldots,n,$
of the independent random variables $\xi_1, \ldots, \xi_n$
satisfy the conditions:
\[
m_{i,k} = {\cal O}(k^{{a_i}k}) \ \mbox{ as }
\ k \to \infty, \ \mbox{ for } \ i=1, \ldots, n,
\]
where $a_1, \ldots, a_n$ are positive constants.
If  $a_1, \ldots, a_n$ are such that
$a_1 + \cdots + a_n \leq 2,$ then the product
$Z_n=\xi_1 \cdots \xi_n$ is M-det. }

\vspace{0.2cm}\noindent
{\bf Proof.} With $m_k={\bf E}[Z_n^k]$ we have,
by the independence of $\xi_i$, that
\[
m_k=m_{1,k} \cdots m_{n,k} ={\cal O}(k^{{a_1}k})
\cdots {\cal O}(k^{{a_n}k}) =  {\cal O}(k^{{a}k})
\ \mbox{ as } \ k \to \infty,
\]
where $a=a_1 + \cdots +a_n$. Since, by assumption, $a \leq 2$, we
apply Theorem 1 to conclude the M-det property of the product
$Z_n.$\hspace{\fill} $\Box$

\vspace{0.2cm}\noindent {\bf Theorem 6.} {\it Suppose that the
growth rates $r_1, \ldots, r_n$ of the moments of each of the random
variables $\xi_1, \ldots, \xi_n$ satisfy
\[
\frac{m_{1,k+1}}{m_{1,k}}={\cal O}((k+1)^{r_1}),
\ \ldots, \ \frac{m_{n,k+1}}{m_{n,k}}={\cal O}((k+1)^{r_n})
\ \mbox{ as } \ k \to \infty,
\]
where $m_{i,k}={\bf E}[\xi_i^k],\ i=1,\ldots,n, \ k=1,2,\ldots. $ If
the rates \ $r_1, \ldots, r_n$ are such that $r_1 + \cdots + r_n
\leq 2,$ then the product \ $Z_n=\xi_1 \cdots \xi_n$ is M-det. }

\vspace{0.2cm}\noindent
{\bf Proof.} Denoting  $m_k={\bf E}[Z_n^k]$ and using the
independence of $\xi_i$, we find that
\[
\frac{m_{k+1}}{m_{k}}=\frac{m_{1,k+1}}{m_{1,k}} \cdots
\frac{m_{n,k+1}}{m_{n,k}} =  {\cal O}((k+1)^{{r}}) \mbox{ as } \ k
\to \infty,
\]
where $r=r_1 + \cdots +r_n$. Since, by assumption, $r \leq 2$, we
refer to Remark 1 and conclude that the product $Z_n$ is
M-det.\hspace{\fill} $\Box$

\vspace{0.2cm}
Let us provide now conditions under which
the product $Z_n$ becomes M-indet.

\vspace{0.2cm}\noindent {\bf Theorem 7.} {\it Consider $n$
independent nonnegative random variables, $\xi_i \sim
F_i,~i=1,2,\ldots, n$, where $n \ge 2.$ Suppose that each $F_i$ is
absolutely continuous with density $f_i>0$ on  $(0,\infty)$  and
that the following conditions are satisfied: \\
(i) At least one (or just one) of the densities
$f_1(x),\ldots,f_n(x)$
is decreasing in $x\ge x_0,$ where $x_0\ge 1$ is a constant.\\
(ii) For each $i=1,2,\dots, n$, there exists  a constant $A_i>0$
such that the density $f_i$ and the tail function $\overline
{F_i}=1-F_i$ satisfy the relation
\begin{eqnarray}
f_i(x)/\overline{F_i}(x)\geq A_i/x~~\hbox{for}~~x\geq x_0,
\end{eqnarray}
and there exist constants $B_i>0,~ \alpha_i>0,$
$\beta_i>0$ and real $\gamma_i$ such that
\begin{eqnarray}
\overline{F_i}(x)\geq B_ix^{\gamma_i}e^{-\alpha_i
x^{\beta_i}}~~ \hbox{for}~~ x\geq x_0.
\end{eqnarray}
If, in addition to (i) and (ii), the parameters $\beta_1, \ldots,
\beta_n$ are such that $\sum_{i=1}^n \frac{1}{\beta_i}>2$, then
the product $Z_n=\xi_1\xi_2\cdots \xi_n$  is M-indet. }\\

To prove Theorem 7, we need the following lemma.

\vspace{0.1cm}\noindent {\bf Lemma 3.} {\it Let $F$ be a
distribution on $\mathbb{R}$ such that (i) it has density $f$ on the
subset $[a, ra]$, where $a>0$ and $r>1$, and (ii) for some constant
$A>0$, ${f(x)}/{\overline{F}(x)}\ge {A}/{x}$ on $[a,ra].$ Then}
$$\int_a^{ra}\frac{f(x)}{x}dx\ge
\left(1-\frac{1}{r}\right)\frac{A}{1+A}\frac{\overline{F}(a)}{a}.$$
\vspace{0.1cm}\noindent {\bf Proof.} By integration by parts, we
have
\begin{eqnarray*}
\int_{a}^{ra}\frac{f(x)}{x}dx&=&-\int_{a}^{ra}\frac{d\overline{F}(x)}{x}
=\frac{\overline{F}(a)}{a}-\frac{\overline{F}(ra)}{ra}-
\int_{a}^{ra}\frac{\overline{F}(x)}{x^2}dx\\
&\ge&\left(1-\frac{1}{r}\right)
\frac{\overline{F}(a)}{a}-\frac{1}{A}\int_{a}^{ra}\frac{f(x)}{x}dx,
\end{eqnarray*}
in which the last inequality is due to the monotonicity of $F$ and
the condition on the failure rate $f/\overline{F}$. Hence the
required conclusion
follows.\hspace{\fill} $\Box$\\

\vspace{0.1cm}\noindent {\bf Proof of Theorem 7.} We may assume,
see condition (i), that $f_n$ is the density which is decreasing
in $x\ge x_0.$ Then, clearly $Z_n$ is nonnegative and its density,
say $h_n,$ can be written as follows: for $x>0$,
\begin{eqnarray*}
&~~& h_n(x)\\
&=&\!\!\!\int_0^{\infty}\!\!\!\int_0^{\infty}\!\!\!\cdots
\!\int_0^{\infty}\frac{f_1(u_1)}{u_1}\frac{f_2(u_2)}{u_2}
\cdots\frac{f_{n-1}(u_{n-1})}{u_{n-1}}f_n\left(\frac{x}{u_1u_2\cdots
u_{n-1}}\right)du_1du_2\cdots du_{n-1}.
\end{eqnarray*}
This representation shows that $h_n(x)>0.$ To prove the M-indet
property of $Z_n$, we will show instead the finiteness of the
Krein quantity of $h_n$. Therefore, we have to estimate the lower
bound of $h_n$.

Since $\sum_{i=1}^n\frac{1}{\beta_i}>2$, we can choose $n-1$
numbers ~$\theta_i\in (0,1),~i=1,2,\ldots, n-1$, such that
$\theta_i<\frac{1}{2\beta_i}$ and
$$1-\frac{1}{2\beta_n} < \sum_{i=1}^{n-1}\theta_i < \min
\left\{1,~\sum_{i=1}^{n-1}\frac{1}{2\beta_i}\right\}.$$ Denote
$\theta_n= 1 - \sum_{i=1}^{n-1}\theta_i,$ define $\theta = \min
\{\theta_1, \ldots, \theta_{n-1}, \theta_n\}$ and let
$x_{\theta}=(2^{n-1}x_0)^{1/\theta}$. Then for each
$x>x_{\theta}$, we take $a_i=x^{\theta_i}\ge x_0$ which implies
\[
x/(2^{n-1}a_1a_2\cdots a_{n-1})=x^{\theta_n}/2^{n-1} \ge x_0.
\]
For these $x$ and $a_i$, we have, by condition (i), the following:
\begin{eqnarray*}
h_n(x)&\geq&\int_{a_1}^{2a_1}\!\!\int_{a_2}^{2a_2}
\!\!\cdots\!\!\int_{a_{n-1}}^{2a_{n-1}}
\frac{f_1(u_1)}{u_1}\frac{f_2(u_2)}{u_2}
\cdots\frac{f_{n-1}(u_{n-1})}{u_{n-1}}\\
&~~&~~~~~~~~~~~~~~~~~~~~\times f_n\left(\frac{x}{u_1u_2\cdots
u_{n-1}}\right)du_1du_2\cdots du_{n-1}\nonumber\\
&\geq&\int_{a_1}^{2a_1}\!\!\int_{a_2}^{2a_2}\!\!\cdots
\int_{a_{n-1}}^{2a_{n-1}}\frac{f_1(u_1)}{u_1}\frac{f_2(u_2)}{u_2}
\cdots\frac{f_{n-1}(u_{n-1})}{u_{n-1}}\\
&~~&~~~~~~~~~~~~~~~~~~~~\times f_n\left(\frac{x}{a_1a_2\cdots
a_{n-1} }\right)du_1du_2\cdots
du_{n-1}\nonumber\\
&=&f_n\left(\frac{x}{a_1a_2\cdots a_{n-1}}\right)\prod_{i=1}^{n-1}
\int_{a_i}^{2a_i}\frac{f_i(u)}{u}du.
\end{eqnarray*}
Then, by Lemma 3 (with $r=2$)
 and (3)--(4),
\begin{eqnarray*}
h_n(x) &\geq&f_n\left(\frac{x}{a_1a_2\cdots a_{n-1}}\right)
\prod_{i=1}^{n-1}\frac{A_i}{2(1+A_i)}\frac{\overline{F_i}(a_i)}{a_i}\\
&\geq& C\left(\frac{x}{a_1a_2\cdots
a_{n-1}}\right)^{\gamma_n-1}\left(\prod_{i=1}^{n-1}a_i^{\gamma_i-1}\right)\\
&~~~&\times \exp\left[{-\sum_{i=1}^{n-1}\alpha_i
a_i^{\beta_i}-\alpha_n\left(\frac{x}{a_1a_2\cdots
a_{n-1}}\right)^{\beta_n}}\right] \\
&=&Cx^{\gamma}\exp\left[-\sum_{i=1}^{n-1}\alpha_i
x^{\theta_i\beta_i}-\alpha_nx^{\theta_n\beta_n}
\right],~~x>x_{\theta},
\end{eqnarray*}
where
$C=2^{1-n}A_n(\prod_{i=1}^{n-1}\frac{A_i}{1+A_i})\prod_{i=1}^nB_i$
and  $\gamma=\sum_{i=1}^{n-1}\theta_i(\gamma_i-1)+
\theta_n(\gamma_{n}-1)$.

 In the exponential factor above we keep
separately the two terms  because of their r\^ole when evaluating
the Krein quantity on $(x_{\theta}, \infty)$ for $h_n.$ Recall
that this is a Stieltjes case and we have the following:
\begin{eqnarray*}
\hbox{K}[h_n]= \int_{x_{\theta}}^{\infty}\frac{-\log
h_n(x^2)}{1+x^2}dx<\infty.
\end{eqnarray*}
 The conclusion about the
finiteness of $\hbox{K}[h_n]$ relies essentially  on the facts
that $2\theta_i\beta_i<1, ~i=1,2,\ldots,n-1$, and
$2\theta_n\beta_n<1.$ Therefore, $Z_n$ is M-indet by Proposition 1
in \cite{p2001}. The proof is complete. \hspace{\fill} $\Box$

\vspace{0.2cm}\noindent
{\bf Example 1.}  For illustration of how to use  Theorem 7,
consider the class of  generalized
gamma distributions. We use the notation $\xi \sim GG(\alpha, \beta,
\gamma)$ if the density function of the random variable $\xi$ is of the form
\begin{eqnarray}
f(x)=cx^{\gamma-1}e^{-\alpha x^{\beta}},\ ~x\geq 0.
\end{eqnarray}
Here $\alpha, \beta, \gamma>0$, $f(0)=0$ if $\gamma\ne 1$, and
$c=\beta\alpha^{\gamma/\beta}/\Gamma(\gamma/\beta)$ is the norming
constant. We have the following statement (see also Theorem 8.4 in
\cite{p2014}).

\vspace{0.2cm}\noindent {\bf Corollary 1.} {\it Suppose  $\xi_1,
\ldots, \xi_n$ are $n$ independent random variables such that
$\xi_i\sim GG(\alpha_i,\beta_i,\gamma_i), \ i=1,\ldots,n$, and let
$Z_n=\xi_1 \cdots \xi_n.$  Then $Z_n$ is M-det if and only if
 $\sum_{i=1}^n\frac{1}{\beta_i}\leq 2$.}

\vspace{0.1cm}\noindent {\bf Proof.} Note that for $\xi\sim
GG(\alpha, \beta, \gamma)$ defined by (5), we have  two
properties, namely: (a) $f(x)/\overline{F}(x)\approx \alpha\beta
x^{\beta-1}$, $\overline{F}(x)\approx
[c/(\alpha\beta)]x^{\gamma-\beta}e^{-\alpha x^{\beta}}$ as
$x\rightarrow \infty$, and (b)
$m_k=\alpha^{-k/\beta}\Gamma((\gamma+k)/\beta))/\Gamma(\gamma/\beta))
={\cal O}(k^{k/\beta})$ as $k\rightarrow \infty$. Hence the
sufficiency part follows from Theorem 1 because ${\bf
E}[Z_n^k]={\cal O}(k^{Bk})$ as $k\to\infty$, where $B=\sum_{i=1}^n
\frac{1}{\beta_i}$, while the necessity part is a consequence of
Theorem 7.\hspace{\fill} $\Box$

\vspace{0.2cm}\noindent {\bf Example 2.} Consider the class of
inverse Gaussian distributions. We say that $X\sim
IG(\mu,\lambda)$ if the density of $X$ is of the form
\begin{eqnarray}
f(x)=\left(\frac{\lambda}{2\pi
x^3}\right)^{1/2}\exp\left[-\frac{\lambda(x-\mu)^2}{2\mu^2x}\right],\
~x>0,
\end{eqnarray}
where $\mu, \lambda>0$ and $f(0)=0.$ If $X\sim IG(\mu,\lambda)$,
then it has a moment generating function. This in turn implies
that the power $Y=X^2$ satisfies Hardy's condition and hence is
M-det. Actually, we have that for real $r$, $X^r$ is M-det if and
only if $|r|\le 2$ (see \cite{s1999}). If $\xi_1$
and $\xi_2$ are two i.i.d.\!\! random variables with density (6),
then the product $Z=\xi_1\xi_2$ is also M-det due to Proposition
1(iii) in \cite{ls2014}. The next result is for products of
non-identically distributed  random variables.

\vspace{0.2cm}\noindent {\bf Corollary 2.} {\it Let $\xi_1\sim
IG(\mu_1,\lambda_1),\ \xi_2\sim IG(\mu_2,\lambda_2)$ and $\eta\sim
Exp(1)$ be
three independent random variables. Then the following statements hold:\\
(i) $Z=\xi_1\eta$ is M-det.\\
(ii) $Z=\xi_1\xi_2$ is M-det.\\
(iii) $Z=\xi_1\xi_2\eta$ is M-indet.}

\vspace{0.1cm}\noindent {\bf Proof.} First, for  $X\sim
F=IG(\mu,\lambda)$, it can be shown (we omit the details) that
 the moment ${\bf E}[X^k]={\cal O}(k^k)$ as $k\to\infty$. Second,
the hazard rate function $r(x)=f(x)/\overline{F}(x)\to
\lambda/(2\mu^2)>0$ as $x\to\infty$. Third,
 the tail function $\overline{F}$ satisfies (4) with the
exponent $\beta=1$. With these three steps we are in a position
to apply Theorems 5 and 7 to confirm the validity of (i) -- (iii) as stated above.
\hspace{\fill} $\Box$

\vspace{0.4cm}\noindent
{\bf 4. Products of Random Variables on} ${\mathbb R} $

\vspace{0.2cm}
We start with two results describing relatively simple
conditions on the
random variables $\xi_1, \ldots, \xi_n$ in order to
guarantee that their product is M-det.
The results are similar to the above Theorems 5 and 6,
however we remember that here we deal with
the Hamburger case, so we work with the even order moments.

\vspace{0.2cm}\noindent
{\bf Theorem 8.} {\it Suppose that the even order moments
$m_{i,2k}={\bf E}[\xi_i^{2k}], \ i=1,\ldots,n,$
of the independent random variables $\xi_1, \ldots, \xi_n$
satisfy the conditions:
\[
m_{i,2k} = {\cal O}((2k)^{{2a_i}k}) \ \mbox{ as } \ k
\to \infty, \ \mbox{ for } \ i=1, \ldots, n,
\]
where $a_1, \ldots, a_n$ are positive constants.
If the parameters $a_1, \ldots, a_n$ are such that
$a_1 + \cdots + a_n \leq 1,$ then the product
$Z_n=\xi_1 \cdots \xi_n$ is M-det. }

\vspace{0.2cm}\noindent
{\bf Proof.} With $m_{2k}={\bf E}[Z_n^{2k}]$ we have,
by the independence of $\xi_i$, that
\[
m_{2k}=m_{1,2k} \cdots m_{n,2k} ={\cal O}((2k)^{{2a_1}k})
\cdots {\cal O}((2k)^{{2a_n}k}) =  {\cal O}((2k)^{{2a}k})
\ \mbox{ as } \ k \to \infty,
\]
where $a=a_1 + \cdots +a_n$. Since, by assumption, $a \leq 1$, we
apply Theorem 3 to  conclude the M-det property of the product
$Z_n.$\hspace{\fill} $\Box$

\vspace{0.2cm}\noindent {\bf Theorem 9.} {\it Suppose  that the
growth rates $r_1, \ldots, r_n$ of the even order moments of each of
the independent random variables $\xi_1, \ldots, \xi_n$ satisfy
\[
\frac{m_{1,2(k+1)}}{m_{1,2k}}={\cal O}((k+1)^{r_1}), \ \ldots,
\ \frac{m_{n,2(k+1)}}{m_{n,2k}}={\cal O}((k+1)^{r_n}) \ \mbox{ as } \ k \to \infty,
\]
where $m_{i,2k}={\bf E}[\xi_i^{2k}], i=1,\ldots,n, \ k=1,2,\ldots. $
If the rates \ $r_1, \ldots, r_n$ are such that $r_1 + \cdots + r_n \leq 2,$
then the product \ $Z_n=\xi_1 \cdots \xi_n$ is M-det. }

\vspace{0.2cm}\noindent
{\bf Proof.} Denoting $m_{2k}={\bf E}[Z_n^{2k}]$ we have,
by the independence of $\xi_i$, that
\begin{eqnarray*}
\frac{m_{2(k+1)}}{m_{2k}}=\frac{m_{1,2(k+1)}}{m_{1,2k}} \cdots
\frac{m_{n,2(k+1)}}{m_{n,2k}} = {\cal O}((k+1)^{r})\ \ \hbox{as}\
k \to \infty,
\end{eqnarray*}
where $r=r_1 + \cdots +r_n$. Since, by assumption, $r \leq 2$,
Theorem 3 together with Remark 4 implies the M-det property of the
product $Z_n.$\hspace{\fill} $\Box$

\vspace{0.2cm} Let us describe now conditions under which the
product $Z_n$ is M-indet.

\vspace{0.2cm}\noindent {\bf Theorem 10.} {\it Consider $n$
independent random variables  $\xi_i \sim F_i, \ i=1,2,\ldots, n$,
 where $n \ge 2$,  and let each $F_i$ be absolutely continuous
 with a symmetric density (about $0$)  $f_i>0$
 on ${\mathbb R}$. Assume further
that the following conditions are satisfied:\\
(i) At least one (or just one) of the densities
$f_1(x),\ldots,f_n(x)$ is decreasing in $x\ge x_0,$ where $x_0\ge 1$ is a constant. \\
(ii) For each $i=1,2,\dots, n$, there exists a constant  $A_i>0$
such that
\begin{eqnarray}
f_i(x)/\overline{F_i}(x)\geq A_i/x~~\hbox{for}~~x\geq x_0,
\end{eqnarray}
and there exist constants $B_i>0,~ \alpha_i>0,$ $\beta_i>0$
and real $\gamma_i$ such that
\begin{eqnarray}
\overline{F_i}(x)\geq B_ix^{\gamma_i}e^{-\alpha_i
x^{\beta_i}}~~ \hbox{for}~~ x\geq x_0.
\end{eqnarray}
If, in addition to the above, $\sum_{i=1}^n\frac{1}{\beta_i}>1$,
then the product $Z_n=\xi_1\xi_2\cdots \xi_n$ is M-indet.}

\vspace{0.1cm}\noindent {\bf Proof.} We may assume, see condition
(i), that $f_n$ is the density which is decreasing in $x\ge x_0.$
Then the density $h_n$ of $Z_n$ is symmetric about 0 (see, e.g.,
\cite{hw1985}) and $h_n$ can be written as follows: \ for $x>0$,
\begin{eqnarray*}
h_n(x)&=&2^{n-1}\int_0^{\infty}\!\!\int_0^{\infty}\!\!\cdots
\int_0^{\infty}\frac{f_1(u_1)}{u_1}\frac{f_2(u_2)}{u_2}
\cdots\frac{f_{n-1}(u_{n-1})}{u_{n-1}}\\
&~~&~~~~~~~~~~~~~~~~~~~~\times f_n\left(\frac{x}{u_1u_2\cdots
u_{n-1}}\right)du_1du_2\cdots du_{n-1}.
\end{eqnarray*}
Hence $h_n(x)>0$. As in the proof of Theorem 7, we will estimate
the lower bound of $h_n$. For completeness we give  detailed
arguments here.

Since  $\sum_{i=1}^n\frac{1}{\beta_i}>1$ by assumption, we can
choose  numbers $\theta_i\in (0,1),~i=1,2,\ldots, n-1$, such that
$\theta_i<\frac{1}{\beta_i}$ and
$$1- \frac{1}{\beta_n}<\sum_{i=1}^{n-1}\theta_i
<\min\left\{1,~\sum_{i=1}^{n-1}\frac{1}{\beta_i}\right\}.$$
Denote
$\theta_n=1-\sum_{i=1}^{n-1}\theta_i$, define
$\theta=\min\{\theta_1,\theta_2,\ldots,\theta_{n-1}, \theta_n\}$
and let $x_{\theta}=(2^{n-1}x_0)^{1/\theta}.$ By taking
$x>x_{\theta}$ and $a_i=x^{\theta_i}$,
 we obtain the following:
\begin{eqnarray*}
h_n(x)&\geq&2^{n-1}\int_{a_1}^{2a_1}\!\!\int_{a_2}^{2a_2}
\!\!\cdots\!\!\int_{a_{n-1}}^{2a_{n-1}}
\frac{f_1(u_1)}{u_1}\frac{f_2(u_2)}{u_2}
\cdots\frac{f_{n-1}(u_{n-1})}{u_{n-1}}\\
&~~&~~~~~~~~~~~~~~~~~~~~\times f_n\left(\frac{x}{u_1u_2\cdots
u_{n-1}}\right)du_1du_2\cdots du_{n-1}\nonumber\\
&\geq&2^{n-1}\int_{a_1}^{2a_1}\!\!\int_{a_2}^{2a_2}\!\!\cdots
\int_{a_{n-1}}^{2a_{n-1}}\frac{f_1(u_1)}{u_1}\frac{f_2(u_2)}{u_2}
\cdots\frac{f_{n-1}(u_{n-1})}{u_{n-1}}\\
&~~&~~~~~~~~~~~~~~~~~~~~\times f_n\left(\frac{x}{a_1a_2\cdots
a_{n-1} }\right)du_1du_2\cdots
du_{n-1}\nonumber\\
&=&2^{n-1}f_n\left(\frac{x}{a_1a_2\cdots
a_{n-1}}\right)\prod_{i=1}^{n-1}
\int_{a_i}^{2a_i}\frac{f_i(u)}{u}du.
\end{eqnarray*}
Applying Lemma 3 and  the above conditions (7) and (8), we further
have
\begin{eqnarray*}
h_n(x) &\geq&2^{n-1}f_n\left(\frac{x}{a_1a_2\cdots a_{n-1}}\right)
\prod_{i=1}^{n-1}\frac{A_i}{2(1+A_i)}\frac{\overline{F_i}(a_i)}{a_i}\\
&\geq& C\left(\frac{x}{a_1a_2\cdots
a_{n-1}}\right)^{\gamma_n-1}\left(\prod_{i=1}^{n-1}a_i^{\gamma_i-1}\right)\\
&~~~&\times \exp\left[{-\sum_{i=1}^{n-1}\alpha_i
a_i^{\beta_i}-\alpha_n\left(\frac{x}{a_1a_2\cdots
a_{n-1}}\right)^{\beta_n}}\right]\\
&=&Cx^{{\tilde \gamma}}\exp\left[-\sum_{i=1}^{n-1}\alpha_i
x^{\theta_i\beta_i}-\alpha_nx^{\theta_n\beta_n}
\right],~~x>x_{\theta}.
\end{eqnarray*}
Here $C=A_n(\prod_{i=1}^{n-1}\frac{A_i}{1+A_i})\prod_{i=1}^nB_i$
and ${\tilde \gamma}=\sum_{i=1}^{n-1}\theta_i(\gamma_i-1)+
\theta_n\,(\gamma_{n}-1)$.

 The next step is to evaluate the Krein
quantity on $|x|>x_{\theta}$ for $h_n$ in this Hamburger case:
\begin{eqnarray*}
\hbox{K}[h_n]= 2\int_{x_{\theta}}^{\infty}\frac{-\log
h_n(x)}{1+x^2}dx<\infty
\end{eqnarray*}
because $\theta_i\beta_i<1$ for each $i=1,2,\ldots, n-1$ and
$\theta_n\,\beta_n<1$.
 Hence, the product $Z_n$ is M-indet by Theorem 2.2 in \cite{p1998}
 (see also \cite{p2001}).  \hspace{\fill} $\Box$

\vspace{0.2cm}\noindent {\bf Example 3.} We now apply Theorem 10
to the product of double generalized gamma random variables. We
write $\xi \sim DGG(\alpha, \beta, \gamma)$ if $\xi$ is a random
variable in $\mathbb R$ with density function of the form
\begin{eqnarray}
f(x)=c|x|^{\gamma-1}e^{-\alpha |x|^{\beta}},\ ~x\in {\mathbb R}.
\end{eqnarray}
Here $\alpha, \beta, \gamma>0$, $f(0)=0$ if $\gamma\ne 1$, and
$c=\beta\alpha^{\gamma/\beta}/(2\Gamma(\gamma/\beta))$ is a norming constant.

\vspace{0.2cm}\noindent {\bf Corollary 3.} {\it Suppose that
$\xi_1, \ldots, \xi_n$ are $n$ independent random variables, and
let $\xi_i\sim DGG(\alpha_i,\beta_i,\gamma_i), \ i=1,2,\ldots,n$.
Then the product $Z_n=\xi_1 \cdots \xi_n$ is M-det if and only if
 $\sum_{i=1}^n\frac{1}{\beta_i}\leq 1$. }

\vspace{0.1cm}\noindent {\bf Proof.} Note that for the moment
$m_{2k}={\bf E}[\xi^{2k}]$ of  $\xi\sim DGG(\alpha, \beta,
\gamma)$, see (9), we have the following relation: $m_{2k}={\cal
O}((2k)^{2k/\beta})$ as $k\rightarrow \infty$. Thus, the
sufficiency part is exactly Theorem 10 in \cite{setal2014}. The
same statement can also be proved by Theorem 8 above. Finally, the
necessity part follows from Theorem 10 (we may redefine $f(0)$ to
be a positive number if necessary).\hspace{\fill} $\Box$

\vspace{0.4cm}\noindent {\bf 5. The Mixed Case}

\vspace{0.2cm} For completeness of our study we need to consider
products of both types of random variables, nonnegative ones and
real ones. Since such a `mixed' product takes values in $\mathbb
R$, this is a Hamburger case, so we can formulate results similar
to Theorems 8 and 9. Since the conditions, the statements and the
arguments are almost as in these two theorems, we do not give
details. Instead, we suggest now a result in which the `mixed'
product $Z_n=\xi_1 \cdots \xi_n$ is M-indet.

\vspace{0.2cm}\noindent
{\bf Theorem 11.} {\it
Given are $n$ independent random variables, such that
 the `first' group, $\xi_1, \ldots, \xi_{n_0}$,
 consists of nonnegative variables, while the variables
in the `second' group,  $\xi_{n_0+1}, \ldots, \xi_n,$ are all with
values in $\mathbb R$, where $1\le n_0<n$. We assume that all
$\xi_i \sim F_i, \ i=1, \ldots, n,$ are absolutely continuous;
denote their densities by $f_i, \ i=1,\ldots,n.$ Assume further
that $f_i(x)>0\ on\ (0,\infty)$, $i=1,\ldots,n_0$, while
$f_j(x)>0\ on\ {\mathbb R},\ j=n_0+1, \ldots, n,$ and are
symmetric.
We require also the following conditions: \\
(i) At least one  of the densities $f_j(x)$, $j=n_0+1,\ldots,n$,
is decreasing in $x\ge x_0,$ where $x_0\ge 1$ is a constant.  \\
(ii) For each $i=1,2,\ldots, n$, there exist  a constant $A_i>0$
such that
\begin{eqnarray*}
f_i(x)/\overline{F_i}(x)\geq A_i/x~~\hbox{for}~~x\geq x_0,
\end{eqnarray*}
and there exist constants $B_i>0,~ \alpha_i>0,$ $\beta_i>0$
and real $\gamma_i$ such that
\begin{eqnarray*}
\overline{F_i}(x)\geq B_ix^{\gamma_i}e^{-\alpha_i
x^{\beta_i}}~~ \hbox{for}~~ x\geq x_0.
\end{eqnarray*}
If, in addition to (i) and (ii), the parameters $\beta_i$ are such
that $\sum_{i=1}^n\frac{1}{\beta_i}>1$, then the product
$Z_n=\xi_1 \cdots \xi_n$ is M-indet.}

\vspace{0.1cm}\noindent {\bf Proof.} The proof is similar to that of
Theorem 10 and can be omitted; the only change is to replace the
coefficient $2^{n-1}$ in the integral form of $h_n$ by
$2^{n-n_0-1}$. \hspace{\fill} $\Box$

\vspace{0.2cm}
As an application of Theorem 11 we derive below two interesting corollaries.

\vspace{0.2cm}\noindent {\bf Corollary 4.} {\it Consider two
independent random variables, $\xi$ and $\eta$, where $\xi \sim
Exp(1)$  and $\eta \sim {\cal N}(0,1)$ (standard normal). Then
$Z=\xi\,\eta$ is M-indet.}

\vspace{0.2cm}\indent In a similar way we arrive also at the
following statement.

\vspace{0.2cm}\noindent {\bf Corollary 5.} {\it (i) The product of two
independent random variables, one chi-square and one normal, is M-indet.\\
(ii) The product of two independent random variables, one inverse Gaussian and one normal,
is M-indet.}

\vspace{0.4cm}\noindent
{\bf Acknowledgements}

\vspace{0.1cm} The work of GDL was partly supported by the
Ministry of Science and Technology of Taiwan (ROC) under Grant NSC
102-2118-M-001-008-MY2.

\vspace{0.1cm}  JS was supported by the International Laboratory
of Quantitative Finance (Higher School of Economics, Moscow) under
Grant No.\,14.12.31.0007 to visit Moscow, March 2014, when a part
of this work was completed. The work was finalized during the
visit of JS to Academia Sinica (Taipei, Taiwan, ROC), May 2014,
with a support from the Emeritus Fellowship provided by Leverhulme
Trust (UK).

\vspace{0.2cm}
\bibliographystyle{plain}

\end{document}